\documentclass[12pt]{amsart}
\usepackage{times,amscd,amssymb,fullpage}
\newtheorem{theorem}{Theorem}
\newtheorem{lemma}[theorem]{Lemma}
\newtheorem{proposition}[theorem]{Proposition}
\newtheorem{corollary}[theorem]{Corollary}
\theoremstyle{definition}
\newtheorem{definition}[theorem]{Definition}

\newtheorem{example}[theorem]{Example}

\theoremstyle{remark}
\newtheorem{remark}[theorem]{Remark}
\newcommand{\FF}{\mathcal{ F}}
\def\varph{{\varphi}}
\def\ra{{\rightarrow}}
\def\lra{{\longrightarrow}}

\def\vGa{{\varGamma}}
\def\bC{{\mathbb C}}
\def\bZ{{\mathbb Z}}
\def\bQ{{\mathbb Q}}
\def\bR{{\mathbb R}}
\def\bH{{\mathbb H}}
\def\inv{{^{-1}}}

\def\M{{\mathcal M}}

\def\TT{{\mathcal T}}

\newcommand\Flux{\operatorname{Flux}}
\newcommand\Symp{\operatorname{Symp}}

\newcommand\Ker{\operatorname{Ker}}

\newcommand\Diff{\operatorname{Diff}}
\newcommand\HEquiv{\operatorname{HEquiv}}

\newcommand\id{\operatorname{Id}}
\catcode`\@=\active 

\begin{document}

\title[Vanishing theorems for flux groups]{Crossed flux homomorphisms and\\ 
vanishing theorems for flux groups}
\author{J.~K\c edra}
\address{Institute of Mathematics, University of Szczecin, Wielkopolska~15, 70-451 Szczecin,
Poland; and Mathematical Sciences, University of Aberdeen, Meston Building, King's College,
Aberdeen AB24 3UE, Scotland, UK}
\email{kedra{\char'100}maths.abdn.ac.uk}
\author{D.~Kotschick}
\address{Mathematisches Institut, Ludwig-Maximilians-Universit\"at M\"unchen,
Theresienstr.~39, 80333 M\"unchen, Germany}
\email{dieter{\char'100}member.ams.org}
\author{S.~Morita}
\address{Department of Mathematical Sciences\\
University of Tokyo \\Komaba, Tokyo 153-8914\\
Japan}
\email{morita{\char'100}ms.u-tokyo.ac.jp}

\date{March 10, 2005, revised August 26, 2005; MSC 2000 classification: primary 57R17, 57R50, 57S05;
secondary 22E65, 53D35}

\begin{abstract}
We study the flux homomorphism for closed forms of arbitrary degree, 
with special emphasis on volume forms and on symplectic forms. The 
volume flux group is an invariant of the underlying manifold, whose 
non-vanishing implies that the manifold resembles one with a circle action 
with homologically essential orbits.
\end{abstract}

\maketitle

\section{Introduction}

\subsection{Flux homomorphisms}
Let $M$ be a closed smooth manifold
and $\alpha$ a closed $p$-form on $M$. We shall denote by $\Diff^{\alpha}$ 
the group of diffeomorphisms of $M$ which preserve $\alpha$, equipped 
with the $C^{\infty}$ topology. Let $\Diff^{\alpha}_{0}$ be its identity component. 
The flux homomorphism associated to $\alpha$ is defined on the universal 
covering of $\Diff^{\alpha}_{0}$ by the formula
\begin{alignat*}{1}
\Flux_{\alpha}\colon\widetilde\Diff^{\alpha}_{0} &\longrightarrow H^{p-1}(M;\bR)\\
\varphi_{t} &\longmapsto \int_{0}^{1}[i_{\dot\varphi_{t}}\alpha] \ dt
\end{alignat*}
for any path $\varphi_{t}$ in $\Diff^{\alpha}_{0}$ with $\varphi_{0}=Id_{M}$.
It is easy to see that the defining integral depends on the path only 
up to homotopy with fixed endpoints.
Identifying an element of the fundamental group of $\Diff^{\alpha}_{0}$ with
the corresponding homotopy class of paths based at the identity in the 
universal covering one obtains a homomorphism 
$$
\Flux_{\alpha}\colon\pi_{1}(\Diff^{\alpha}_{0}) \longrightarrow H^{p-1}(M;\bR)
$$
whose image is the flux group $\Gamma_{\alpha}$ associated with $\alpha$.
The flux homomorphism descends to a homomorphism defined on $\Diff^{\alpha}_{0}$, 
also called the flux:
$$
\Flux_{\alpha}\colon\Diff^{\alpha}_{0}\longrightarrow 
H^{p-1}(M;\bR)/\Gamma_{\alpha} \ .
$$

There are a number of general questions one can ask about this situation, such as 
whether the flux group $\Gamma_{\alpha}$ is trivial, or at least discrete, and whether 
the flux homomorphism can be extended from $\Diff^{\alpha}_{0}$ to the whole group 
$\Diff^{\alpha}$. As far as we know, these questions have only been considered in the 
literature in the case when $\alpha$ is a symplectic form, see for 
example~\cite{Banyaga,B,Calabi,Kedra,KM,LMPflux,LMP,LO,McDuff,MS} and the papers cited there.

It is the aim of this paper to discuss these questions in some generality. 
In Section~\ref{s:general} we shall show how certain arguments used in~\cite{KM} for the case of 
symplectic forms, mostly on surfaces, can be adapted to the general case, proving 
the following result:
\begin{theorem}\label{t:mainA}
    Suppose that the total space of every $M$-bundle with structure group $\Diff^{\alpha}$ has a 
    cohomology class restricting to $[\alpha]$ on the fiber. Then the flux homomorphism
    $$
    \Flux_{\alpha}\colon \widetilde{\Diff^{\alpha}_{0}}\longrightarrow H^{p-1}(M;\bR)
    $$ 
    vanishes on $\pi_{1}(\Diff^{\alpha}_{0})$ and extends to a crossed homomorphism 
    $$
    \widetilde{\Flux_{\alpha}}\colon \Diff^{\alpha}\longrightarrow H^{p-1}(M;\bR) \ .
    $$
    \end{theorem}
The crossed flux homomorphism $\widetilde{\Flux_{\alpha}}$ is a cocycle representing a cohomology 
class with coefficients in $H^{p-1}(M;\bR)$ on the group $\Diff^{\alpha}$ considered as a discrete group, 
which extends the cohomology class on $\Diff^{\alpha}_0$ given by the flux homomorphism. 

There are many other situations in which we only prove the vanishing of the flux group
$\Gamma_{\alpha}$, without exhibiting a crossed homomorphism extending the flux homomorphism.
An important instance of this occurs when $\alpha$ represents a bounded cohomology class, 
see Theorem~\ref{t:bounded}.

\subsection{Volume flux}
In Sections~\ref{s:volume} and~\ref{s:entropy} we shall study the flux homomorphism for volume forms $\mu$.
We begin by showing that a smooth circle action gives rise to a non-trivial volume flux group if and only if
its orbits are homologically essential in real homology. In fact, these are the only known examples of non-trivial
volume flux, and it might be possible that there are no others. In dimensions $1$ and $2$
the only closed manifolds with non-trivial volume flux groups are $S^1$, respectively $T^2$. For these manifolds
$\Diff^{\mu}_0$ is homotopy equivalent to the manifold itself, and the loops in $\Diff^{\mu}_0$ with non-trivial
flux are generated by smooth circle actions. Modulo the Poincar\'e conjecture, this last statement is also true in 
dimension $3$, as we will show in Section~\ref{3-folds} by proving the following:
\begin{theorem}\label{t:SeifertA}
    Let $M$ be a closed oriented $3$-manifold without any fake cells.
    If $M$ has non-zero volume flux group $\Gamma_{\mu}$, then $M$ is a 
    Seifert fiber space and a multiple of every loop 
    with non-zero volume flux is realized by a fixed-point-free circle action on $M$. 
    \end{theorem}

Our guiding paradigm then is to extend to manifolds with a non-trivial volume flux group the known restrictions on 
manifolds with homologically essential circle actions. The results we have described so far show that a non-trivial
volume flux group implies the vanishing of all real characteristic numbers, see Corollary~\ref{t:volume}, and of the
simplicial volume, see Corollary~\ref{t:simplicial}. In the case of fixed-point-free circle actions, these  results are 
consequences of the vanishing of the minimal volume. Recall that the minimal volume, introduced by Gromov 
in~\cite{Gromov}, is defined by 
$$
\textrm{MinVol}(M)=\inf\{ Vol(M,g) \ \vert \ g\in Met(M) \ 
\textrm{with} \ \vert K_{g}\vert\leq 1 \} \ , 
$$
where $K_{g}$ denotes the sectional curvature of $g$. It is known that this is a very sensitive invariant of $M$, 
which depends on the smooth structure in an essential way. Even the vanishing or non-vanishing of the  minimal 
volume depends subtly on the smooth structure~\cite{entropies}. In the presence of a fixed-point-free smooth circle 
action on $M$ the minimal volume vanishes because one can shrink a suitable invariant metric in the direction of 
the orbits, and thus collapse $M$ while keeping the sectional curvatures bounded. It is tempting to speculate that a 
non-trivial volume flux group is enough to imply the vanishing of the minimal volume, and this would follow if one could
prove that circle actions also account for non-trivial volume flux groups in dimensions $>3$.
In order to describe further results in this direction, we need to recall certain notions of entropy.

\subsection{Flux groups and entropies} 
Let $M$ be a connected closed oriented manifold of dimension $n$. Elaborating on ideas of Gromov~\cite{Gromov},
the following lower bounds for the minimal volume of $M$ have been proved, compare~\cite{BCG0,entropies,PP}:
\begin{equation}\label{e:main}
    \frac{n^{n/2}}{n!} \vert\vert M \vert\vert \leq 2^{n}n^{n/2} T(M) \leq 
    \lambda(M)^{n}\leq h(M)^{n}\leq (n-1)^{n}\text{MinVol}(M) \ .
    \end{equation}
Here $T(M)$ is the spherical volume introduced by Besson, Courtois and Gallot~\cite{BCG0}, and $h(M)$ is the 
minimal topological entropy of geodesic flows on $M$. We will have nothing to say about these two invariants,
but shall be concerned with the other three quantities in~\eqref{e:main}, namely the simplicial volume $\vert\vert M \vert\vert$, 
the minimal volume entropy or asymptotic volume $\lambda (M)$, and the minimal volume $\text{MinVol}(M)$.

For a Riemannian metric $g$ on $M$ consider the lift $\tilde g$ to the 
universal covering $\tilde M$. For an arbitrary basepoint $p\in\tilde M$ 
consider the limit 
$$
\lambda (M,g) = \lim_{R\to\infty}\frac{log Vol(B(p,R))}{R} \ ,
$$
where $B(p,R)$ is the ball of radius $R$ around $p$ in $\tilde M$ with respect 
to $\tilde g$, and the volume is taken with respect to $\tilde g$ as well. After 
earlier work by Efremovich, Shvarts, Milnor~\cite{Mi} and others, Manning 
showed that the limit exists and is independent of $p$. 
It follows from~\cite{Mi} that $\lambda (M,g)>0$ if and only if $\pi_{1}(M)$ has exponential growth. 
We call $\lambda (M,g)$ the volume entropy of the metric $g$, and define the 
minimal volume entropy or asymptotic volume of $M$ to be 
$$
\lambda (M) = \inf\{\lambda (M,g)\ \vert \  g\in Met(M) \ 
\textrm{with} \ Vol(M,g)=1 \} \ .
$$
This sometimes vanishes even when $\lambda (M,g)>0$ for every $g$.
The normalization of the total volume is necessary because of the 
scaling properties of $\lambda (M,g)$.
Babenko proved that the minimal volume entropy $\lambda (M)$ is invariant
under homotopy equivalences, and also under certain bordisms over $B\pi_1(M)$, 
see~\cite{B1,B2}. 

Of all these invariants, only the simplicial volume is known to be multiplicative in 
coverings. As it will be convenient to allow ourselves passage to finite coverings,
we make the following definition in the spirit of~\cite{MT}, compare also~\cite{BCG0}.

\begin{definition}
Let $I$ be an invariant of $n$-dimensional closed manifolds. Then define
$$
I^*(M)=\inf\left\{ \frac{I(N)}{d} \ \vert \ N \ \textrm{a} \ d-\textrm{sheeted covering of}  \ M \right\} \ , 
$$
where the infimum is taken over all finite coverings of $M$.
\end{definition}
Clearly $I^*(M)\leq I(M)$, and if $I(M)\leq J(M)$ for all $M$, then $I^*(M)\leq J^*(M)$. Thus,
\eqref{e:main} implies
\begin{equation}\label{e:main2}
    \frac{n^{n/2}}{n!} \vert\vert M \vert\vert \leq 2^{n}n^{n/2} T^*(M) \leq 
    \lambda^*(M)^{n}\leq h^*(M)^{n}\leq (n-1)^{n}\text{MinVol}^*(M) \ .
    \end{equation}

As we mentioned already, the existence of a smooth circle action without fixed points on $M$ implies that its minimal
volume vanishes. Therefore, all the quantities in~\eqref{e:main} and~\eqref{e:main2}
vanish. In the case of a non-trivial volume flux group, we are not able to prove the vanishing of the minimal
volume, but, in Section~\ref{s:entropy} we shall prove the following weaker result:
\begin{theorem}\label{t:entropyA}
    Let $M$ be a closed oriented manifold with non-vanishing volume flux 
    group $\Gamma_{\mu}$. Then $M$ has a finite covering $\bar M$ 
    whose volume entropy $\lambda (\bar M)$ vanishes. In particular, $\lambda^*(M)=0$.
\end{theorem}
The proof uses in an essential way the bordism invariance of the minimal volume entropy $\lambda (M)$,
proved by Babenko~\cite{B2}.

\subsection{Symplectic flux}
In Section~\ref{s:powers} we shall consider symplectic forms and their powers, for which we obtain 
generalizations of some results previously proved in~\cite{McDuff,Kedra,KM}. 
Unlike in the case of volume forms, where the flux group is always discrete for 
purely topological reasons, the discreteness of the symplectic flux group was 
an open problem until very recently. This issue, first raised by Banyaga, has
just been resolved by Ono's proof~\cite{Ono} using methods of hard symplectic 
topology. Our results in Section~\ref{s:powers}, like those of some of the references
mentioned above, show that very often the symplectic flux group actually vanishes.

\section{General properties of the flux and crossed flux homomorphisms}\label{s:general}

For cyclic isotopies representing elements of the fundamental group of 
$\Diff^{\alpha}_{0}$ we can reformulate the definition of the flux as follows.
\begin{lemma}\label{l:tori}
If $\varphi_{t}$ is a closed loop representing an element of 
$\pi_{1}(\Diff^{\alpha}_{0})$ and $c$ is a cycle representing a 
homology class in $H_{p-1}(M,\bZ)$, then up to sign we have
$$
\langle \Flux_{\alpha}(\varphi_{t}),[c]\rangle \ = \langle 
[\alpha],[\varphi_{t}(c)]\rangle \ ,
$$
where $\varphi_{t}(c)$ denotes the cycle swept out by the loop of 
diffeomorphisms $\varphi_{t}$ applied to $c$.
\end{lemma}
\begin{proof}
    This is immediate from the definition of the flux; compare~\cite{McDuff}.
    \end{proof}
As every real cohomology class can be detected by mapping closed 
oriented manifolds $\Sigma$ into $M$, this lemma shows that the flux 
group of $\alpha$ is detected by the evaluation of the pullback of 
$[\alpha]$ on products $S^{1}\times\Sigma$. More precisely, define 
$\phi\colon S^{1}\times\Sigma\rightarrow M$ by 
$\phi(t,x)=\varphi_{t}(f(x))$, where $f\colon\Sigma\rightarrow M$ 
is a representative for (a multiple of) $[c]$. Then 
$$
\langle \Flux_{\alpha}(\varphi_{t}),[\Sigma]\rangle \ = \langle 
\phi^{*}[\alpha],[S^{1}\times\Sigma]\rangle \ .
$$
One should not be misled by this discussion into thinking that the 
flux group depends only on the cohomology class of $[\alpha]$, 
because which loops in $\Diff_{0}$ can be deformed to essential loops 
in $\Diff^{\alpha}_{0}$ depends on $\alpha$ itself, and not just on 
its cohomology class.

\begin{remark}\label{r:ss}
Let $[\varphi_{t}]\in\pi_1(\Diff^{\alpha}_{0})$, and denote by 
$E_{\varphi_{t}}\to S^2$ the $M$-bundle corresponding
to $\varphi_{t}$ via the clutching construction.
It is known~\cite{LO,Oprea} that the flux homomorphism 
$\Flux_{\alpha}\colon\pi_1(\Diff^{\alpha}_{0})\to H^{p-1}(M,\bR)$
is equal to the evaluation on $[\alpha]$ of the differential 
$$
\partial_{\varphi_{t}}^*\colon H^{*}(M;\bR)\longrightarrow H^{*-1}(M;\bR)
$$
in the cohomology spectral sequence associated to $E_{\varphi_{t}}$. 
More precisely,
$$
\left <\Flux_{\alpha}(\varphi_{t}),c\right > = 
\left <\partial_{\varphi_{t}}^*[\alpha], c \right >
$$
for all $c\in H_{p-1}(M)$.
\end{remark}

The diffeomorphism group $\Diff^{\alpha}$ acts by conjugation on itself and 
on the universal covering of its identity component. It also acts on 
cohomology, and this latter action factors through the mapping class 
group $\M_{\alpha}$ with respect to $\alpha$, defined to be the quotient group 
$\Diff^{\alpha}/\Diff^{\alpha}_{0}$. Our first observation is that the 
flux is equivariant with respect to these actions:
\begin{lemma}\label{l:equiv}
The flux homomorphism $\Flux_{\alpha}\colon\widetilde\Diff^{\alpha}_{0} 
\longrightarrow H^{p-1}(M;\bR)$ is
equivariant with respect to the natural actions of $\Diff^{\alpha}$.
In other words, for any two elements
$\psi\in\Diff^{\alpha}$ and $\varphi_{t}\in\widetilde\Diff^{\alpha}_{0}$, 
we have the identity
$$
\Flux(\psi\varphi_{t}\psi\inv)=\bar{\psi}(\Flux(\varphi_{t}))
$$
where $\bar{\psi}\in\M_{\alpha}$ denotes the mapping class of $\psi$ and
$\M_{\alpha}$ acts on $H^{p-1}(M;\bR)$ from the left by the rule
$\bar{\psi}(w)=(\bar{\psi}^{-1})^*(w)$ for $w\in H^{p-1}(M;\bR)$.
\end{lemma}
This follows immediatly from the definition of the flux and the chain rule.  

The lemma suggests that one should not expect an extension of the flux 
to $\Diff^{\alpha}$ to exist as a homomorphism, but rather as a 
crossed homomorphism with respect to this action of $\Diff^{\alpha}$ 
on cohomology. Indeed, in certain situations we shall prove the 
existence of an extension as a crossed homomorphism.

Consider the extension
\begin{equation*}\label{eqn:mom}
1\lra\Diff^{\alpha}_{0}\lra\Diff^{\alpha}
\lra\M_{\alpha}\lra 1 
\end{equation*}
and its associated exact sequence of cohomology groups of discrete 
groups:
\begin{align*}
0 & \lra H^1(\M_{\alpha};H^{p-1}(M;\bR)/\vGa_{\alpha})  \lra
H^1(\Diff^{\alpha};H^{p-1}(M;\bR)/\vGa_{\alpha})\\
& \lra 
H^1(\Diff^{\alpha}_{0};H^{p-1}(M;\bR)/\vGa_{\alpha})^{\M_{\alpha}}
\overset{\delta}{\lra}
H^2(\M_{\alpha};H^{p-1}(M;\bR)/\vGa_{\alpha})\lra
\end{align*}
Lemma~\ref{l:equiv} shows that we can think of the flux homomorphism 
as an element
$$
\Flux_{\alpha}\in 
H^1(\Diff^{\alpha}_{0};H^{p-1}(M;\bR)/\vGa_{\alpha})^{\M_{\alpha}} \ .
$$

Extending the flux to $\Diff^{\alpha}$ as a crossed homomorphism is 
equivalent to the vanishing of $\delta(\Flux_{\alpha})$ in the above exact 
sequence. We now examine this issue in detail. 

For a foliated $M$-bundle $E\rightarrow B$ whose total holonomy is 
contained in $\Diff^{\alpha}$ we have a transverse invariant class 
$a\in H^{p}(E;\bR)$ defined as follows. Pulling back $\alpha$ from $M$ to 
$\tilde B\times M$, we obtain an invariant form $\tilde\alpha$ which descends 
to $E=(\tilde B\times M)/\pi_{1}(B)$ as a closed form. We denote its 
cohomology class by $a$.
\begin{lemma}\label{l:a}
Let $I=[0,1]$. For any $\varph\in\Diff^{\alpha}_{0}$ let
$\pi\colon M_{\varph}\ra S^1$ be the foliated $M$-bundle
over $S^1$ with monodromy $\varph$. It is the quotient
of $M\times I$ by the equivalence relation
$(p,0)\sim (\varph(p),1)$. By assumption,
there is an isotopy $\varph_{t}\in \Diff^{\alpha}_{0}$
such that $\varph_{0}=\id$ and $\varph_{1}=\varph$.
Let $f\colon M_{\varph}\rightarrow M\times S^1$ be the induced
diffeomorphism given by the correspondence
$$
M_{\varph}\ni (p,t) \longmapsto (\varph_{t}^{-1}(p),t)\in M\times
S^1 \ .
$$

Then the transverse invariant class $a\in H^{p}(M_{\varph};\bR)$ 
is equal to
\begin{align*}
[\alpha] + \Flux_{\alpha}(\varph_{t})\otimes\nu\in &H^{p}(M\times S^1;\bR)\\
\cong &H^{p}(M;\bR)\oplus (H^{p-1}(M;\bR)\otimes H^1(S^1;\bR))
\end{align*}
under the above isomorphism, where $\nu\in H^1(S^1;\bR)$
denotes the fundamental cohomology class of $S^1$.
\end{lemma}
\begin{proof}
The horizontal foliation on $M_{\varph}$ is induced from the trivial
foliation on $M\times I$. Hence the transverse invariant class $a$ 
is represented by the form $p_{1}^* \alpha$ on $M\times I$, where 
$p_{1}\colon M\times I\ra M$ denotes the projection to the first factor. 
It is clear that the $H^{p}(M;\bR)$-component of $a$ is equal to 
$[\alpha]$, so that we only need to prove that for any $(p-1)$-cycle 
$c \subset M$, the value of $a$ on the cycle 
$f\inv(c\times S^1) \subset M_{\varph}$ is equal to
$\Flux(\varph_{t})([c])$ where $[c]\in H_{p-1}(M;\bZ)$ denotes
the homology class of $c$. Now on $M\times I$, the above
cycle is expressed as the image of the map
$$
c\times I\ni (q,t)\longmapsto (\varph_t(q),t)\in M\times I
$$
because $f\inv(q,t)=(\varph_t(q),t)\ ((q,t)\in M\times S^1)$.
Hence the required value is equal to the integral of $\alpha$ over
the image of the map
$$
c\times I\ni (q,t)\longmapsto \varph_t(q)\in M \ .
$$
But this is exactly equal to the value of $\Flux(\varph_{t})$ on the
homology class represented by the cycle $c\subset M$.
This completes the proof.
\end{proof}

For the formulation of our result about extensions of the flux 
homomorphism as a crossed homomorphism we use the following 
definition:
\begin{definition}\label{d:ext}
    Let $M$ be a closed manifold and $G\subset\Diff (M)$ a subgroup. 
    We say that a cohomology class $c\in H^{*}(M)$ extends 
    $G$-universally if there is a class $b$ on the total space 
    of any $M$-bundle with structure group $G$ restricting to $c$ on 
    the fibers.
    \end{definition}
With this terminology we have the following precise version of Theorem~\ref{t:mainA}
from the introduction:
\begin{theorem}\label{t:main}
    Let $G\subset\Diff^{\alpha}$ be a subgroup. If $[\alpha]$ 
    extends $G$-universally, then the flux homomorphism
    $$
    \Flux_{\alpha}\colon \widetilde{G_{0}}\longrightarrow H^{p-1}(M;\bR)
    $$ 
    vanishes on $\pi_{1}(G_{0})$ and extends to a crossed 
    homomorphism 
    $$
    \widetilde{\Flux_{\alpha}}\colon G\longrightarrow H^{p-1}(M;\bR) \ .
    $$
    \end{theorem}
\begin{proof}
By assumption $[\alpha]$ extends $G$-universally, so in particular it 
extends to the $M$-bundle over $S^{2}$ given by the clutching 
construction for a loop in $G_{0}$. Therefore, $\Flux_{\alpha}$ 
vanishes on $\pi_{1}(G_{0})$, cf.~Remark~\ref{r:ss}.
    
Let $BG^{\delta}$ be the classifying space of $G$ 
considered as a discrete group, and denote by 
$$
\pi\colon EG^{\delta}\lra BG^{\delta}
$$
the universal foliated $M$-bundle over $BG^{\delta}$
with total holonomy group in $G$. Let $b\in H^{p}(EG^{\delta};\bR)$ 
be a universal extension of $[\alpha]$, and consider the difference
$$
u=a- b \in H^{p}(EG^{\delta};\bR) \ ,
$$
where $a$ denotes the transverse invariant class represented by
the global $p$-form $\tilde\alpha$ on $EG^{\delta}$ which
restricts to $\alpha$ on each fiber.
The restriction of $u$ to the fiber vanishes, so that,
in the spectral sequence $\{E^{p,q}_r\}$
for the real cohomology, we have the natural projection
\begin{align*}
P\colon\Ker&\left(H^{p}(EG^{\delta};\bR)\ra
H^{p}(M;\bR)\right) \ni u\\
&\lra P(u)\in E_{\infty}^{1,p-1}\subset E_{2}^{1,p-1}= 
H^{1}(BG^{\delta};H^{p-1}(M;\bR)) \ .
\end{align*}
Now Lemma~\ref{l:a} implies that the restriction of $P(u)$ to the 
identity component of $G$ coincides with the flux homomorphism:
\begin{equation*}
P(u)=\Flux_{\alpha}\colon G_{0}\lra H^{p-1}(M;\bR) \ .
\end{equation*}
Thus we see that $P(u)$ defines an extension of the flux homomorphism 
as a cohomology class whose representing cocycles are crossed 
homomorphisms.
\end{proof}

A particular case where the class $[\alpha]$ does extend universally 
is when it represents some characteristic class for $M$:
\begin{corollary}\label{c:main}
    Suppose $[\alpha]\in H^{p}(M;\bR)$ is a non-zero multiple of a 
    polynomial in the Euler and Pontryagin classes of $M$. Then the 
    flux group $\Gamma_{\alpha}$ vanishes and the flux homomorphism 
    extends as a crossed homomorphism
    $$
    \widetilde{\Flux_{\alpha}}\colon\Diff^{\alpha}\longrightarrow H^{p-1}(M;\bR) \ .
    $$
    \end{corollary}
\begin{proof}
    For any $M$-bundle $E\rightarrow B$ consider the tangent bundle 
    along the fibers. Its characteristic classes extend the 
    characteristic classes of $TM$ from the fiber to the total space.
    \end{proof}

If $\alpha$ defines a geometric structure on $M$, then $\Diff^{\alpha}$ 
acts by automorphisms of this structure and preserves its characteristic 
classes. We shall consider the case of a symplectic structure in 
Section~\ref{s:powers} below. Another instance of this is the case of foliations:
\begin{example}
    Suppose $\alpha$ is of constant rank, and $T\FF\subset TM$ is its 
    kernel. Then $\Diff^{\alpha}$ preserves $T\FF$, and its 
    characteristic classes extend $\Diff^{\alpha}$-universally. 
    Therefore, the flux group vanishes and the flux homomorphism 
    extends as a crossed homomorphism if $[\alpha]$ is a non-zero 
    multiple of a polynomial in the Euler and Pontryagin classes of 
    $T\FF$ and of $TM/T\FF$.
    \end{example}

In general one cannot expect that the extension of the flux homomorphism to
a crossed homomorphism is unique (if it exists at all). However, it was proved in~\cite{KM}
that for the case of a symplectic form, equivalently an area form, on a surface of
genus $\geq 2$, the extension is unique.

There is another mechanism which can force the vanishing of the flux group $\Gamma_{\alpha}$, 
stemming from Gromov's notion of bounded cohomology~\cite{Gromov}. As usual, we say that a 
real cohomology class $[\alpha]\in H^p(M;\bR)$ is a bounded class, if it has a representative which 
is bounded as a functional on the set of singular simplices. This means that the class in the image 
of the comparison map from bounded to usual cohomology:
$$
H^p_b(M;\bR)\longrightarrow H^p(M;\bR) \ .
$$
    
\begin{theorem}\label{t:bounded}
Suppose the closed $p$-form $\alpha$ represents a bounded cohomology class. 
Then the flux group $\Gamma_{\alpha}$ vanishes.
\end{theorem}
\begin{proof}
Suppose the flux group $\Gamma_{\alpha}$ does not vanish. Then according
to Lemma~\ref{l:tori} there is a smooth map $\phi\colon S^1\times\Sigma^{p-1}\rightarrow M$
for which $\phi^*\alpha$ has non-zero integral over $S^1\times\Sigma$. As $[\alpha]$
is assumed to be a bounded class, so is $\phi^*[\alpha]$. This means that the cohomology
generator in top degree on $S^1\times\Sigma$ is bounded, equivalently the simplicial
volume of $S^1\times\Sigma$ is non-zero. This is clearly false, as $S^1\times\Sigma$
maps to itself with arbitrarily large degree.
\end{proof}
In this case we obtain the vanishing of the flux group although we do not have a crossed
homomorphism extending the flux homomorphism.
    
\section{Volume-preserving diffeomorphisms}\label{s:volume}

In this section we consider the flux of $\alpha$ when $\alpha = \mu$ is 
a volume form on $M$. In this case Moser's celebrated result~\cite{Moser} 
implies that $\Diff^{\mu}$ is weakly homotopy equivalent to the full 
diffeomorphism group $\Diff (M)$ of $M$. Moreover, applying 
Lemma~\ref{l:tori} to a homotopy of volume forms of equal total 
volume, we deduce that, up to normalization, $\Gamma_{\mu}$ is 
independent of $\mu$. Furthermore, in this situation the mapping 
class group $\M_{\mu}$ does not depend on $\mu$, as it equals 
$\Diff (M)/\Diff_{0} (M)$.

To see some examples and get a feel for what results to expect, we 
first consider loops of diffeomorphisms generated by circle actions.

\subsection{Circle actions}
Suppose we are given a smooth effective circle action on an oriented closed manifold 
$M$. By averaging we can always construct an invariant volume form 
$\mu$, so that we have a non-trivial homomorphism $S^{1}\rightarrow 
\Diff^{\mu}_{0}$. We can easily characterize the non-triviality of the 
flux on the image of this homomorphism:
\begin{proposition}\label{p:circle}
    A circle action gives rise to a nonzero element in $\Gamma_{\mu}$ 
    if and only if its orbits are nonzero in real homology.
    \end{proposition}
    \begin{proof}
	Let $X$ be the vector field generating the $S^{1}$-action. Then 
$L_{X}\mu = 0$, and $i_{X}\mu$ is a closed $S^{1}$-invariant $(n-1)$-form 
representing the volume flux evaluated on the image of $\pi_{1}(S^{1})$. 
This form is also a defining form for the (singular) foliation defined by 
the orbits. 

Choose closed $S^{1}$-invariant $1$-forms $\alpha_{1},\ldots,\alpha_{k}$ 
representing a basis for the first de Rham cohomology. Then the wedge 
products $\alpha_{i}\wedge i_{X}\mu$ are $S^{1}$-invariant as well, 
and their cohomology classes vanish for all $i\in\{1,\ldots,k\}$ if 
and only if the flux vanishes in cohomology. But $\alpha_{i}\wedge 
i_{X}\mu$ is constant along each orbit, and vanishes if and only if 
$\alpha_{i}(X)$ vanishes along the orbit. As all the orbits are homologous 
to each other and fill out the manifold, the flux can only vanish if 
$\alpha_{i}(X)$ vanishes identically for all $i$, which is equivalent to 
the orbits being null-homologous.
	(Note that it is not possible that $\alpha_{i}\wedge i_{X}\mu$ is 
	exact but not identically zero.)
	\end{proof}
	
If the orbits are non-trivial in homology, then the action has no fixed points. 
A closed one-form $\alpha_{i}$ representing an integral class and evaluating 
non-trivially on an orbit has $\alpha_{i}(X)\neq 0$ everywhere, and therefore 
defines a smooth fibration over $S^{1}$ with fibers transverse to the circle action.
The finiteness of the isotropy groups of the circle action implies finiteness of the 
monodromy of the fibration over $S^{1}$. Thus our manifold has a finite cover 
which splits off $S^{1}$ smoothly and equivariantly, with the standard circle action. 
This gives a differential-geometric proof of the Conner--Raymond 
theorem~\cite{CR}, originally proved by topological means.
(Compare~\cite{FM} for similar arguments.)

Given any fixed-point-free circle action on $M$, Gromov~\cite{Gromov} 
showed that the minimal volume of $M$ vanishes. {\it A fortiori}, the 
simplicial volume and the real characteristic numbers of $M$ vanish. 
Specializing the general theorems of the previous section to volume
forms allows us to extend these vanishing results from circle actions to 
non-trivial volume flux groups. Further, it is well-known that the orbits of 
circle actions represent central elements of the fundamental group acting 
trivially on homotopy groups, see for example Browder--Hsiang~\cite{BH} or 
Appendix~2 in~\cite{Gromov}. We shall generalize these statements in 
Theorem~\ref{p:vol-flux}.

In the case of $3$-manifolds there are more precise results. Namely, it 
was proved by Epstein that any fixed-point-free circle action on a closed 
$3$-manifold occurs by rotating the fibers of a Seifert fibration. Moreover, 
if the orbits are non-trivial in real homology, it is easy to see that the 
Euler class of the fibration is trivial, so that $M$ is finitely covered by a 
product of a surface with the circle. We shall show in Theorem~\ref{t:Seifert} 
below that these circle fibrations account for all non-trivial volume flux 
elements on $3$-manifolds.

\subsection{Topological consequences of non-vanishing volume flux groups}
We begin with a characterization of the non-triviality of the volume 
flux group, together with some homotopical constraints.
\begin{theorem}\label{p:vol-flux}
Let $M$ be any closed $n$-manifold with volume form $\mu$.
\begin{enumerate}
\item The volume flux group $\Gamma_{\mu}$ is trivial if and only if
the evaluation map $ev\colon\Diff^{\mu}_0 \to M$ induces the
trivial map on the first real homology.
\item If $\Gamma_{\mu}\neq 0$, then $ev_{*}\colon\pi_{1}(\Diff^{\mu}_0) \to 
\pi_{1}(M)$ has an infinite image, which acts trivially on the 
homotopy groups of $M$. In particular, the center of $\pi_{1}(M)$ is 
infinite.
\end{enumerate}
\end{theorem}
\begin{proof}
Suppose $\varphi_{t} \in \pi_1(\Diff^{\mu}_0)$ 
and denote by $M\to E_{\varphi_{t}}\to S^2$ the bundle associated
to $\varphi_{t}$ by the clutching construction.
According to Remark~\ref{r:ss}, the non-vanishing of 
$\Flux_{\mu}(\varphi_{t})$ is equivalent to the non-vanishing of
$\partial^*_{\varphi_{t}}(\mu )$. Applying Poincar\'e
duality in $E_{\varphi_{t}}$, this in turn is equivalent to the 
non-triviality of the differential
$\partial_{\varphi_{t}}\colon H_{0}(M)\to H_{1}(M)$.
On the other hand, $\partial_{\varphi_{t}}[pt] = ev_*(\varphi_{t})$,
where $ev_*\colon H_1(\Diff^{\mu}_0)\to H_1(M)$ is the map
induced by the evaluation at the point $pt\in M$.
This proves the first claim.

Considering the evaluation on the fundamental group, we conclude from 
what we proved above that it has infinite image. It is a general 
property of the image of the evaluation that it acts trivially on all 
homotopy groups. This was first noticed by Gottlieb~\cite{G0},
compare also Theorem~2.2 in~\cite{Oprea}. 
\end{proof}

In the case of volume forms, Corollary~\ref{c:main} gives the following:
\begin{corollary}\label{t:volume}
    Let $M$ be a closed oriented manifold of dimension $2n$, and $\mu$ 
    a volume form on $M$. If $M$ has a nonzero real characteristic number, 
    then the flux group $\Gamma_{\mu}$ is trivial, and the flux homomorphism 
    $\Flux_{\mu}$ extends to a crossed homomorphism 
    $$
    \widetilde\Flux_{\mu}\colon\Diff^{\mu}\longrightarrow 
    H^{2n-1}(M;\bR) \ .
    $$
    \end{corollary}
It was proved in Theorem~2 of~\cite{KM} that the cohomology class of the 
extension $\widetilde\Flux_{\mu}$ is uniquely determined in 
$H^1(\Diff^{\mu};H^{1}(M;\bR))$ if $\mu$ is a volume form on a surface 
of genus $g\geq 2$. When $M$ has dimension at least $4$, and has two different 
nonzero characteristic numbers, for example the Euler characteristic and the 
signature, then it may happen that these two choices give rise to different extensions 
of the flux associated with a volume form. In cohomology, the 
difference between any two such extensions is in 
$H^1(\M_{M};H^{2n-1}(M;\bR))$, where $\M_{M}=\Diff(M)/\Diff_{0}(M)$ 
is the mapping class group of $M$. 

Corollary~\ref{t:volume} only applies to even-dimensional manifolds. 
In all dimensions, we have the following special case of Theorem~\ref{t:bounded}:
\begin{corollary}\label{t:simplicial}
    Let $M$ be a closed oriented manifold with nonzero simplicial 
    volume. Then the volume flux group $\Gamma_{\mu}$ is trivial for 
    every volume form $\mu$.
    \end{corollary}
    
\subsection{The case of $3$-manifolds}\label{3-folds}
We can now prove a precise version of Theorem~\ref{t:SeifertA} mentioned in the 
introduction:
\begin{theorem}\label{t:Seifert}
    Let $M$ be a closed oriented $3$-manifold without any fake cells.
    If $M$ has nonzero volume flux group $\Gamma_{\mu}$, then $M$ is a 
    Seifert fiber space. In particular, its minimal volume vanishes.
    Moreover, a multiple of every loop $\varphi_{t}\in\pi_{1}(\Diff^{\mu}_0)$ 
    with nonzero volume flux is realized by a fixed-point-free circle action on $M$. 
    \end{theorem}
\begin{proof}
    As $\pi_{1}(M)$ has non-trivial center by Theorem~\ref{p:vol-flux}, it 
    is indecomposable as a free product. Therefore, in the Kneser--Milnor 
    prime decomposition~\cite{M} of $M$, all summands but one are simply connected. As 
    $M$ contains no fake cells by assumption, we conclude that it is prime. Thus, 
    either $M$ is $S^{1}\times S^{2}$, which is a Seifert fibration in the 
    obvious way, or it is irreducible.
    
    If $M$ is irreducible, then because its first Betti number is 
    positive, it is also sufficiently large\footnote{Thus $M$ 
    is ``Haken''.}, meaning that $M$ contains an incompressible 
    surface, cf.~\cite{Jaco} p.~35. Now it is a theorem of 
    Waldhausen~\cite{W} that a closed irreducible sufficiently large 
    $3$-manifold $M$ such that $\pi_{1}(M)$ has non-trivial center is 
    Seifert fibered.
    
    By shrinking a suitable invariant metric in the direction of the circle 
    action, one sees that the minimal volume vanishes, 
    cf.~\cite{Gromov}.
    
    On a Seifert manifold every element of the center of the 
    fundamental group is, up to a multiple, represented by the fiber 
    of a Seifert fibering, cf.~\cite{Jaco} p.~92/93 or~\cite{Scott}, and 
    therefore by a circle action. Thus a multiple of the evaluation of a loop 
    $\varphi_{t}\in\pi_{1}(\Diff^{\mu}_0)$ with non-trivial volume flux is homotopic
    in $M$ to the orbits of a circle action. It remains to show that $\varphi_{t}$ 
    and the loop given by the circle action are homotopic in $\Diff^{\mu}_0$, before 
    we apply the evaluation. For this we distinguish the two cases we encountered
    above: either $M$ is $S^{1}\times S^{2}$, or it is irreducible.
    
    If $M$ is $S^{1}\times S^{2}$, we can use the work of Hatcher~\cite{H1,H2},
    who determined the homotopy type of $\Diff (S^{1}\times S^{2})$ completely. This shows in particular
    that the free part of the fundamental group of $\Diff_0$, and therefore of 
    $\Diff^{\mu}_0$, is generated by the circle action $S^1\hookrightarrow\Diff_0(S^{1}\times S^{2})$
    given by rotation of the first factor.
    
    If $M$ is irreducible, consider the space $\HEquiv (M)$ of self-homotopy equivalences 
    of $M$, endowed with the compact-open topology, and the following composition of 
    maps:
    $$
    \Diff^{\mu}_0(M)\stackrel{i}{\longrightarrow}\Diff_0(M)\stackrel{j}{\longrightarrow}\HEquiv_0(M)
    \stackrel{ev}{\longrightarrow} M \ ,
    $$ 
    where $\HEquiv_0(M)$ denotes the connected component of the identity. The first map is a 
    weak homotopy equivalence by Moser's theorem~\cite{Moser}. For an irreducible sufficiently 
    large $3$-manifold $M$, Laudenbach~\cite{L} proved that $j_*$ is an isomorphism on fundamental
    groups\footnote{Laudenbach~\cite{L} assumed the Smale conjecture, subsequently proved 
    by Hatcher~\cite{H2}.}. An irreducible $3$-manifold with infinite fundamental group is 
    aspherical by the sphere theorem, see~\cite{M}, and it is true for any aspherical manifold that 
    $ev\colon\HEquiv_0(M)\longrightarrow M$ induces an isomorphism between $\pi_1(\HEquiv_0(M))$ 
    and the center of $\pi_1(M)$, compare~\cite{G0}. Thus, two loops in $\Diff^{\mu}_0$ having 
    homotopic evaluations in $M$ are homotopic in $\Diff^{\mu}_0$. This completes the proof of the theorem.
    \end{proof}
\begin{remark}
The Seifert fibered $3$-manifolds occuring in the theorem carry Thurston geometries of type 
$S^{2}\times\bR$ in the case of $S^1\times S^2$, and of type $\bR^{3}$ or $\bH^{2}\times\bR$
in the irreducible case, compare~\cite{Scott}.  The center of $\pi_1(M)$ is $\bZ^3$ if $M$ is $T^3$,
and is $\bZ$ otherwise. For Seifert manifolds which are not circle bundles over a surface, the 
generator of the center can not be realised by a circle action, so that passing to multiples is 
unavoidable. We shall see a similar phenomenon in higher dimensions in Example~\ref{ex:7D} below.
\end{remark}

\subsection{Back to higher dimensions}
The proof of Theorem~\ref{t:Seifert} has a partial generalization to 
higher dimensions:
\begin{theorem}\label{t:irred}
    Let $M$ be a closed $n$-manifold with $\Gamma_{\mu}\neq 0$. If $M$ 
    is homotopy equivalent to a connected sum $M_{1}\# M_{2}$ then one 
    of the $M_{i}$ is a homotopy sphere.
    \end{theorem}
    \begin{proof}
	 As $\pi_{1}(M)$ has non-trivial center by Theorem~\ref{p:vol-flux}, it 
    is indecomposable as a free product. Therefore, one of the $M_{i}$, 
    say $M_{1}$, is simply connected. In particular, 
    $H_{1}(M_{1};\bZ)=0$. If $M_{1}$ is not a homotopy sphere, then 
    there is a smallest $k\leq n-2$ for which $H_{k}(M_{1};\bZ)$ does not 
    vanish. By the Hurewicz theorem $\pi_{k}(M_{1})\cong 
    H_{k}(M_{1};\bZ)\neq 0$.
    
    Now the universal cover of $M$ is obtained from the universal cover of 
    $M_{2}$ by connected summing with infinitely many copies of 
    $M_{1}$. Every non-trivial element of 
    $\pi_{1}(M)\cong\pi_{1}(M_{2})$ acts non-trivially on $H_{k}(\tilde 
    M)$ by permuting the different summands coming from the different 
    copies of $M_{1}$. This shows that every element of the 
    fundamental group acts non-trivially on $\pi_{k}(M)$, contradicting 
    the second part of Theorem~\ref{p:vol-flux}.
	\end{proof}
\begin{remark}
Theorem~\ref{t:irred} is new even in the case when the non-triviality of the 
volume flux group arises from a circle action. See~\cite{FM} for partial results
in this direction. 
\end{remark}
The following example shows that Theorem~\ref{t:irred} is sharp.
\begin{example}\label{ex:7D}
    Let $M = (S^{1}\times S^{6})\#\Sigma^{7}$, with $\Sigma$ a 
    homotopy $7$-sphere. Now every $\Sigma$ is a 
    twisted sphere, i.~e.~it is of the form $D^{7}\bigcup_{\psi}D^{7}$ 
    for some $\psi$ in the mapping class group of $S^{6}=\partial D^{7}$. 
    But then $M$ is just the mapping torus of $\psi$, and as the mapping 
    class group of $S^{6}$ is finite (of order $28$), we conclude that $M$ 
    fibers over $S^{1}$ with finite monodromy ($=\psi$). There is a 
    fixed-point-free circle action transverse to this fibration 
    generating a non-trivial volume flux group by 
    Proposition~\ref{p:circle}.
    
    Note that the generator of $\pi_{1}(M)\cong\bZ$ cannot always be 
    realized by the orbits of a circle action on $M$, so that passing to 
    multiples is unavoidable. Indeed, if an $S^{1}$-action on $M$  
    surjects $\pi_{1}(S^{1})$ onto $\pi_{1}(M)$, then all the orbits 
    have trivial stabilizer, because their homotopy classes are primitive. 
    Then the quotient map $M\rightarrow M/S^{1}$ is a smooth circle bundle 
    over a homotopy $6$-sphere. This bundle is trivial, and as there are 
    no exotic $6$-spheres we conclude that $M$ is diffeomorphic to 
    $S^{1}\times S^{6}$. But according to Browder~\cite{Browder}, 
    Corollary~2.8, we can choose $\Sigma$ in such a way that $M$ is 
    not diffeomorphic to $S^{1}\times S^{6}$. 
    \end{example}
    This example shows that the topological manifold $S^{1}\times S^{6}$ 
    has several distinct smooth structures, all of which admit 
    fixed-point-free circle actions. There are also pairs of homeomorphic 
    manifolds for which one has a free smooth circle action, and the other
    one has no smooth circle action at all:
    \begin{example}\label{ex:tori}
    Let $M = T^7\#\Sigma^{7}$, with $\Sigma$ a homotopy $7$-sphere. 
    Whenever $\Sigma$ is not the standard sphere, Assadi and Burghelea~\cite{AB}
    showed that $M$ admits no effective smooth circle action. 
        \end{example}
    We do not know whether the volume flux group is non-trivial in this case, or not.
See Section~\ref{ss} below for further discussion.

\section{Entropy and volume flux}\label{s:entropy}

We have seen that a non-trivial volume flux group forces the vanishing of 
the simplicial volume, and the vanishing of all real characteristic numbers.
In the case when the volume flux comes from a smooth circle action,
these vanishing results are consequences of the vanishing of the minimal
volume, compare~\eqref{e:main}. One might therefore speculate that the non-vanishing of the 
volume flux may imply the vanishing of the minimal volume. As we are
not able to prove this, we consider the intermediate invariants from~\eqref{e:main},
which interpolate between the simplicial volume and the minimal volume.
The strongest vanishing result we can prove about them is Theorem~\ref{t:entropyA}.
The remainder of this section is occupied by the proof of this theorem.

If the volume flux group $\Gamma_{\mu}$ for $(M,\mu)$ is non-trivial, then 
    by Theorem~\ref{p:vol-flux} the evaluation at a point of the corresponding 
    loop $\varphi_t$ in $\Diff^{\mu}$ gives us a loop which is of infinite order 
    in the center of $\pi_{1}(M)$ and in $H_{1}(M;\bZ)$. After replacing $M$ by a  
    finite cover, we may assume that this element is primitive in $H_{1}(M;\bZ)/tor$, 
    so that the fundamental group of $M$ splits as a direct product 
    $\pi_{1}(M)\cong\bZ\times\pi$ with the generator of the first factor corresponding 
    to the evaluation of our loop of diffeomorphisms, compare~\cite{O}. In this situation, 
    Gottlieb~\cite{Gottlieb} and independently Oprea~\cite{O} proved a homotopical 
    analogue of the Conner--Raymond splitting theorem, showing that $M$ is 
    homotopy equivalent to $S^{1}\times Y$, where $Y$ is the homotopy fiber of a 
    map $f\colon M\rightarrow S^{1}$ inducing the projection 
    $p_{1}\colon\pi_{1}(M)\rightarrow\bZ$ onto the first factor of the fundamental group. 
    If the homotopy type $Y$ can be represented by a closed oriented $(n-1)$-manifold, 
    then we conclude the proof of Theorem~\ref{t:entropyA} by noting that $\lambda (S^{1}\times Y)$ 
    vanishes because of the obvious circle action, and the minimal entropy is known to be 
    homotopy-invariant by a result of Babenko~\cite{B1}.
    
    Regardless what the homotopy fiber $Y$ is, we can proceed as follows.
    Choose a smooth map $f\colon M\rightarrow S^{1}$ with $f_{*}=p_{1}$ and let 
    $F'$ be a regular fiber. Then $\pi_{1} (F')$ surjects onto $\pi\cong \Ker f_{*}$, 
    and we can modify $F'$ by ambient surgery\footnote{At this point it is useful 
    to assume that $dim (F')>2$, equivalently $dim (M)\geq 4$. This is no loss of 
    generality, as we have a much stronger conclusion for small dimensions by a 
    different argument, see Theorem~\ref{t:Seifert}.} inside $M$ to obtain an 
    embedded submanifold $F\subset M$ in the same homology class, such that 
    $\pi_{1}(F)\cong\pi$, and the inclusion induces an isomorphism between 
    $\pi_{1}(F)$ and $0\times\pi\subset\bZ\times\pi\cong\pi_{1}(M)$. 
    Consider then the map $\Phi$ given by the composition
    \begin{alignat*}{3}
    S^{1} \times F &\stackrel{i}{\longrightarrow} S^1\times M &&\stackrel{Id\times\varphi}{\longrightarrow} S^1\times M &&\stackrel{\pi_2}{\longrightarrow} M\\
    (t,x) &\longmapsto (t,x) &&\longmapsto (t,\varphi_{t}(x)) &&\longmapsto\varphi_{t}(x)  \ ,
    \end{alignat*}
    where the first map is the inclusion, and the composition of the second and third maps is the evaluation.
\begin{lemma}\label{l:class}
The map $\Phi$ has the following properties:
\begin{enumerate}
\item It induces an isomorphism on fundamental groups.
\item It has degree one.
\item It pulls back the tangent bundle of $M$ to the tangent bundle of $S^{1}\times F$.
\end{enumerate}
\end{lemma}
\begin{proof}
The first claim is clear from the construction of $F$ and $\Phi$. The second claim 
follows from the fact that $F\subset M$ has algebraic intersection number $=1$ with 
the evaluation loop.

For the third claim, consider the factorization $\Phi = \pi_2\circ(Id\times\varphi) \circ i$. 
The diffeomorphism $Id\times\varphi$ pulls back $\pi_2^*TM$ to itself. But this bundle
restricts to the image of $i$ as $\bR\oplus TF$, which proves the claim. 
\end{proof}

    Let $c\colon M\rightarrow B\pi_{1}(M)$ be the classifying map for the universal 
    cover of $M$, and consider the classes $[M,c]$ and $[S^{1}\times F,c\circ\Phi ]$ 
    in the bordism group $\Omega_{n}(B\pi_{1}(M))$. If these bordism classes 
    agree, then there is a bordism $[W',\alpha']$ between them. It follows that 
    $\alpha_{*}'\colon\pi_{1}(W')\rightarrow\pi_{1}(M)$ is surjective, and we can 
    modify the bordism by surgery in the interior of $W'$ so as to obtain a new 
    bordism $[W,\alpha]$ for which $\alpha_{*}\colon\pi_{1}(W)\rightarrow\pi_{1}(M)$ 
    is an isomorphism. This new bordism has the property that the inclusion of each 
    boundary component into $W$ induces an isomorphism on fundamental groups,
    which is the defining property of an $R$-cobordism in the terminology of Babenko~\cite{B2}.
    The result of~\cite{B2} is that $R$-cobordant manifolds have the same minimal 
    volume entropy. As before, the volume entropy of $S^{1}\times F$ vanishes 
    because its minimal volume vanishes courtesy of the circle action, cf.~\eqref{e:main}.
    
    It remains to prove that the bordism classes $[M,c]$ and $[S^{1}\times F,c\circ\Phi ]$ 
    agree---or to deal with their failure to do so. 
    Consider first the case when the integral homology of $B\pi_1(M)$ has no odd-order
    torsion. Then the bordism spectral sequence for $B\pi_1(M)$ is trivial, and we have
    an isomorphism
    \begin{equation}\label{e:bordism}
    \Omega_{n}(B\pi_{1}(M))\cong\bigoplus_{i=0}^n H_{i}(B\pi_{1}(M);\Omega_{n-i}(\star)) \ ,
    \end{equation}
    compare Theorem~15.2 in~\cite{CF}, or~\cite{Stong}.
    The elements of the summands on the right-hand side are detected by the collection of all Pontryagin 
    and Stiefel-Whitney numbers twisted by cohomology classes on $B\pi_{1}(M)$.
    Lemma~\ref{l:class} shows that these twisted characteristic numbers agree
    for $[M,c]$ and $[S^{1}\times F,c\circ\Phi ]$, so these two bordism classes in 
    $\Omega_{n}(B\pi_{1}(M))$ agree.
    
   In the general case, when the homology of $B\pi_1(M)$ is allowed to have odd-order
    torsion, we can still find the required bordism between $[M,c]$ and $[S^{1}\times F,c\circ\Phi ]$
    after passing to a suitable finite cover induced from a finite cover of $S^1$. For this we only have to prove
    that on such a cover the map $\Phi\colon S^1\times F\longrightarrow M$ is bordant to the 
    identity of $M$.
\begin{proposition}\label{p:surgery}
In the above situation $M$ has a finite covering induced from a finite covering of $S^1$ via the map
$M\longrightarrow S^1$, such that the lift of $\Phi$ to the corresponding covering of $S^1\times F$
by itself is bordant to the identity of the target.
\end{proposition}
\begin{proof}
We shall use the language of surgery theory.
First we recall a few basic facts from this theory, see e.g.~\cite{MM}.

Let $O=\lim_{n\to\infty} O_n$ denote the
infinite orthogonal group and let $\mathrm{B}O$ denote its
classifying space. Also let $G_n$ denote the group consisting of
homotopy equivalences of $S^{n-1}$ equipped with the compact-open
topology. Let $\mathrm{B}G$ be the classifying space of
$G=\lim_{n\to\infty} G_n$, which classifies stable
isomorphism classes of spherical fibrations.
The homotopy groups $\pi_i (G)$ can be canonically identified
with the stable homotopy groups of spheres
$$
\pi_i (G)\cong \lim_{k\to\infty} \pi_{i+k} (S^k) \ .
$$
Hence $\pi_i (G)$ is finite for all $i$.
There exists a canonical map
$\mathrm{B}O\to \mathrm{B}G$, which corresponds to associating
the unit sphere bundle to a vector bundle.
We have a fibration
$$
G/O\lra \mathrm{B}O\lra \mathrm{B}G
$$
which can be extended to the left as
$$
G\lra G/O\lra \mathrm{B}O\lra \mathrm{B}G \ .
$$
Now let $X$ be a closed oriented smooth manifold and let $\nu(X)$ be
its stable normal bundle. A commutative diagram
\begin{equation}
\begin{CD}
\nu(N) @>{\tilde f}>> \nu(X)\\
@VVV @VVV\\
N @>{f}>> X,
\end{CD}
\label{eqn:bdle}
\end{equation}
where $f\colon N\to X$ is a degree $1$ map from a closed
oriented manifold $N$ to $X$ and $\tilde f\colon\nu(N)\to\nu(X)$
is a bundle map, is called a {\it normal map} over $X$.
There exists a notion of {\it normal cobordism}, which is an
equivalence relation on the set of all normal maps over $X$.
The set of all the normal cobordism classes of normal maps
over $X$ is denoted by $NM_O(X)$. The map
$$
\sigma\colon NM_O(X)\lra [X,G/O]
$$
is called the surgery map. It is known to be a bijection by the Pontryagin--Thom
construction, see~\cite{MM}, Theorem 2.23. A normal map
$(\tilde f, f)$ as in~\eqref{eqn:bdle} is normally cobordant
to the identity map of $X$ if and only if $\sigma$
sends the corresponding class $[(\tilde f, f)]\in NM_O(X)$
to the homotopy class of the constant map in $[X,G/O]$.

Now we apply this general framework to our
map $\Phi\colon S^1\times F \to M$.
First of all, we already know from Lemma~\ref{l:class} that $\Phi$ is a
tangential equivalence, so that it is covered by a bundle
map $\tilde\Phi\colon T(S^1\times F)\ra TM$. Taking the stable
normal bundles instead of the tangent bundles, we obtain
a normal map
\begin{equation}
\begin{CD}
\nu(S^1\times F) @>{\tilde \Phi^{\nu}}>> \nu(M)\\
@VVV @VVV\\
S^1\times F @>{\Phi}>> M
\end{CD}
\label{eqn:ff}
\end{equation}
over $M$.
Consider the surgery map
$$
\sigma\colon NM_0(M)\lra [M, G/O]
$$
and let $\alpha\in NM_0(M)$ be the element represented
by the normal map~\eqref{eqn:ff}.
Our task is to show that, after passing to a suitable
finite cover, the image $\sigma(\alpha)$ is the
homotopy class in $[M, G/O]$ represented by
the constant map.

Passing to finite coverings along $S^1$ means that we consider
$\Phi_k\colon S^1 \times F \longrightarrow M_k$, where the domain is
the product of $F$ with the standard $k$-fold cover of the circle, and $M_k$
is the corresponding covering of $M$ pulled back from the circle.
Since $\Phi_k$ has degree $1$, we see that
$$
\Phi_k^*\colon H^i(M_k) \longrightarrow H^i(S^1 \times F)
$$
is injective with any coefficients. Moreover, by restriction we have an injection
$$ 
\Phi_k^*\colon H^i(M_k,F) \longrightarrow H^i(S^1 \times F, pt. \times F)\cong \iota\otimes H^{i-1}(F)  \ , 
$$
where $\iota$ denotes the generator of $H^1(S^1)$ and $F$ sits inside $M_k$ in the obvious way
for any $k$.




Now $\Phi$ is a tangential equivalence, and therefore a map $M\longrightarrow G/O$ 
representing $\sigma(\alpha)$ can be lifted to $g\colon M\longrightarrow G$.
Moreover, this map is constant on $F \subset M$ because $\Phi$ is the identity on $F$. 
It suffices to show that $g$ is homotopic to the constant map, at least after we pass
to a suitable finite covering. The obstructions to $g$ being null-homotopic are contained in
$$
H^i(M,F;\pi_i(G)) \subset \iota\otimes H^{i-1}(F; \pi_i(G)) \ .
$$
But the group $\pi_i(G)$ is finite for every $i$, and if we pass to the finite covering
$M_k$, then 
the classifying map of the new surgery problem is just the composition
$$
M_k \longrightarrow M \stackrel{g}{\longrightarrow} G \ ,
$$
as can be seen by inspecting the construction of the surgery map, cf.~pp. 42-43 of~\cite{MM}.
Therefore we can kill the obstructions by taking suitable finite covers along $S^1$.
In detail, the coverings send $\iota$ to $k \iota$, and since the obstructions lie in
$\iota\otimes H^{i-1}(F)$ with finite coefficients, we can kill all the obstructions. This means that
the normal map  $S^1 \times F \longrightarrow M_k$ is normally bordant to the identity
of $M_k$.
\end{proof}
The proof of Proposition~\ref{p:surgery} shows that $S^1 \times F$ and $M_k$ are bordant 
over $B\pi_1 (M_k)$.
Together with the preceding discussion, this completes the proof of Theorem~\ref{t:entropyA}.

\begin{remark}
    An alternative approach to the general case proceeds by observing that~\eqref{e:bordism}
    always holds after tensoring with $\bQ$, cf.~\cite{CF,Stong}. As $[M,c]$ and 
    $[S^{1}\times F,c\circ\Phi ]$ have the same twisted Pontryagin numbers, their difference 
    is rationally zero-bordant. This means that for some $k>0$ there is a bordism between 
    the $k$-fold connected sums $kM$ and $k(S^{1}\times F)$ endowed with the 
    corresponding maps to $B\pi_{1}(M)$. These maps induce the diagonal map
    $D\colon\pi_1(M)\star\ldots\star\pi_1(M)\longrightarrow\pi_1(M)$ on fundamental 
    groups. Unfortunately it is unclear whether this can be arranged to be an $R$-cobordism in the sense of 
    Babenko~\cite{B2}. If this is possible, then a slightly different proof of Theorem~\ref{t:entropyA} can be given 
    as follows. Babenko's theorem~\cite{B2} implies that the minimal asymptotic exponential 
    volume growth rates of the covers of $kM$ and of $k(S^{1}\times F)$ with fundamental 
    groups $Ker(D)$ are equal. (These are not the minimal volume entropies, because these 
    covers are not the universal covers.)
    Now by a result of Paternain and Petean~\cite{PP}, Theorem 5.9, the circle 
    action on $S^{1}\times F$ gives rise to a $\TT$-structure on the connected 
    sum $k(S^{1}\times F)$. Another result of the same authors, Theorem A 
    in~\cite{PP}, shows that the minimal topological entropy $h$ vanishes for 
    any manifold with a $\TT$-structure. {\it A fortiori}, the minimal volume 
    entropy of $k(S^{1}\times F)$ vanishes, compare~\eqref{e:main}. This 
    implies that the intermediate cover of $k(S^{1}\times F)$ with fundamental
    group $Ker(D)$ also has slow volume growth. By the above discussion
    we have this conclusion also for the cover of $kM$ with fundamental 
    group $Ker(D)$. This cover essentially 
    contains a copy of the universal cover of $M$, which therefore has 
    small minimal asymptotic exponential volume growth rate.
    Thus, $\lambda(M)=0$.
\end{remark}

\section{Powers of a symplectic form}\label{s:powers}

In this section we consider the case when $\alpha$ is a power 
$\omega^{k}$ of a symplectic form $\omega$ on $M$, with $M$ of 
dimension $2n$. It is clear that $\Diff^{\alpha}$ contains the 
symplectomorphism group $\Symp = \Diff^{\omega}$, but is usually strictly 
larger when $k>1$. In order to obtain a result parallel to 
Corollary~\ref{t:volume}, we want to use the Chern classes of the tangent 
bundle along the fibers in the universal foliated $M$-bundle. This means 
that instead of $\Diff^{\alpha}$ we should only consider a smaller group which 
preserves the homotopy class of an almost complex structure compatible 
with $\omega$. We will simply take the symplectomorphism group and 
consider the $k$-flux
$$
\Flux_{k}\colon\widetilde\Symp_{0}\longrightarrow H^{2k-1}(M;\bR) \ ,
$$
which is the restriction of the flux with respect to $\omega^{k}$ to 
the symplectomorphism group. We denote by $\Gamma_{k}$ the image of 
$\pi_{1}(\Symp_{0})$ under the $k$-flux.

The groups $\Gamma_{k}$ for different values of $k$ are related to 
each other by the equation
\begin{equation}\label{e:kflux}
        \Flux_{k}(\varphi_{t}) = k \ \Flux_{1}(\varphi_{t})\wedge\omega^{k-1} \ ,
\end{equation}
        which is immediate from the definition of the flux and the identity 
        $i_{X}(\omega^{k}) = k \ i_{X}\omega\wedge\omega^{k-1}$. Thus 
        $\Gamma_{k}$ is the image of of the usual symplectic flux 
        group $\Gamma_{1}=\Gamma_{\omega}$ under multiplication by 
        $k\omega^{k-1}$. Taking $k=n$, we can use this to draw 
        consequences about the symplectic flux group from our results 
        about the volume flux. Note that $\Gamma_{n}$ is not really 
        the volume flux group, because we are only considering 
        $\Symp$, and not the usually larger $\Diff^{\omega^{n}}$. Nevertheless, the same 
        arguments apply to prove the following:
\begin{theorem}\label{t:kfl}
    Let $(M,\omega)$ be a closed symplectic manifold of dimension $2n$ that 
    satisfies one of the following conditions:
    \begin{enumerate}
    \item the evaluation map $ev\colon\Symp_0\to M$ induces the trivial 
    map on the first real homology, or
    \item the fundamental group $\pi_{1}(M)$ has finite center, or
    \item $M$ has a nonzero real characteristic number, or has 
    nonzero renormalized minimal volume entropy $\lambda^* (M)$, or
    \item $M$ is homotopy equivalent to a connected sum in which 
    neither summand is a homotopy sphere. 
    \end{enumerate}
    Then the symplectic flux group $\Gamma_{\omega}\subset 
    H^{1}(M;\bR)$ is in the kernel of the multiplication map 
    $$
    [\omega]^{n-1}\colon H^{1}(M;\bR)\rightarrow H^{2n-1}(M;\bR) \ .
    $$
    
    In particular, if $(M,\omega)$ also satisfies the hard Lefschetz 
    condition, then the symplectic flux group vanishes.
    \end{theorem}

Instead of taking the maximal power of the symplectic form, we can 
consider smaller powers, and we can also use the Chern classes of an 
almost complex structure compatible with the symplectic form.
Theorem~\ref{t:main} has the following immediate consequence:
\begin{corollary}\label{t:symp}
    Suppose that $[\omega^{k}]\in H^{2k}(M;\bR)$ is proportional to a 
    polynomial in the Chern classes of $(M,\omega)$. Then
    $\Gamma_{k}=0$, and $\Flux_{k}$ extends to a crossed homomorphism
    $$
    \widetilde\Flux_{k}\colon\Symp\longrightarrow H^{2k-1}(M;\bR) \ .
    $$
    \end{corollary}
This result has an antecedent in McDuff's paper~\cite{McDuff}. The case $k=1$ 
was proved in~\cite{KM}.

As before, using~\eqref{e:kflux} we obtain the following consequence. 
The case when $k=n$ was previously proved in~\cite{Kedra}.
\begin{corollary}
    Suppose that $[\omega^{k}]\in H^{2k}(M;\bR)$ is proportional to a 
    polynomial in the Chern classes of $(M,\omega)$. Then the symplectic flux 
    group $\Gamma_{1}=\Gamma_{\omega}\subset H^{1}(M;\bR)$ 
    is in the kernel of the multiplication map 
    $$
    [\omega]^{k-1}\colon H^{1}(M;\bR)\longrightarrow 
    H^{2k-1}(M;\bR) \ .
    $$
    If $(M,\omega)$ satisfies a weak form of the Lefschetz property, 
    namely if multiplication by $[\omega]^{k-1}$ is injective, 
    then the usual symplectic flux group is trivial.
    \end{corollary}

\begin{example}
    Consider $M=F\times S^{2}$, where $F$ is a surface of genus $g\neq 
    1$. Then every cohomology class with nonzero square in 
    $H^{2}(M;\bR)$ is realised by a split symplectic form, with the 
    symplectic area of the factors scaled suitably. For all these 
    symplectic forms the Chern classes are the same, namely 
    $c_{1}=(2-2g)P.D.[S^{2}]+2P.D.[F]$ and $c_{2}=(4-4g)P.D.[M]$. For 
    those symplectic forms $\omega$ whose cohomology class is a 
    multiple of $c_{1}$, the case $k=1$ of Corollary~\ref{t:symp} 
    implies the triviality of the flux group $\Gamma_{\omega}$. 
    When $[\omega]$ is not a multiple of $c_{1}$, we can use the case 
    $k=2$ and the fact that $c_{1}^{2}$ and $c_{2}$ are nonzero to 
    conclude that $\Gamma_{2}$ is trivial\footnote{Alternatively we 
    can use the fact that $\pi_{1}(M)$ has trivial center.}. 
    As $M$ satisfies the hard Lefschetz property for every $\omega$, 
    we again conclude the vanishing of $\Gamma_{\omega}$.
\end{example}

\begin{example}
    Let $M=F\times S^{2}$ as before, with $g\neq 1$. Then the 
    non-vanishing of the Chern numbers $c_{1}^{2}$ and $c_{2}$ gives 
    rise to two potentially different extensions
    $$
    \widetilde\Flux_{2}\colon\Symp\longrightarrow H^{3}(M;\bR) \ .
    $$
    The difference between them corresponds to the difference 
    $c_{1}^{2}(\xi)-2c_{2}(\xi)\in H^{4}(E\Symp;\bR)$, which 
    restricts trivially to the fiber $M$. There are symplectic bundles with 
    fiber $M$ which show that this difference class is non-trivial if 
    $g=g(F)\geq 3$. Namely, let $X\rightarrow B$ be an $F$-bundle 
    with nonzero signature.
    Then $X\times S^{2}$ is an $M$-bundle over $B$ for 
    which $c_{1}^{2}(\xi)-2c_{2}(\xi)\neq 0\in H^{4}(X\times S^{2};\bR)$.
    
    However, the two extensions of $\Flux_{2}$ are essentially the same.
    As $M$ satisfies the hard Lefschetz property, these extensions of 
    $\Flux_{2}$ are given by extensions of the usual flux homomorphism
    multiplied by the symplectic form. But the extensions of the usual flux
    homomorphism here come from $F$, where we know that the extension
    is unique as a cohomology class, see~\cite{KM}.
    \end{example}

\begin{example}
    Consider now $M=T^{2}\times S^{2}$ with a split symplectic form. Then 
    clearly the flux group is non-trivial. However, if we pass from $M$ to 
    its blowup $\hat M = M\#\overline{\bC P^{2}}$, then $\hat M$ is reducible 
    (and has nonzero Chern numbers), so that $\Gamma_{2}$ must vanish. As 
    $\hat M$ satisfies the hard Lefschetz property, we conclude that the 
    usual flux group $\Gamma_{1}=\Gamma_{\omega}$ also vanishes.
    \end{example}
    
    Theorem~\ref{t:bounded} has the following consequence for the 
    symplectomorphism groups.
\begin{corollary}\label{t:sb}
    Suppose that $[\omega]^{k}$ is a bounded cohomology class. Then 
    $\Gamma_{k}=0$. In particular, if $\omega$ represents a bounded class, 
    then the usual symplectic flux group is trivial.
    \end{corollary}

It is interesting to compare the above vanishing results for the symplectic flux group with the 
following:
\begin{proposition}\label{p:P}{\rm (cf.~\cite{LMPflux,Polterovich})}
Assume that $(M,\omega)$ is symplectically aspherical, i.~e. $\omega\vert_{\pi_2(M)}=0$,
and that $\pi_1(M)$ has finite center. Then $\Gamma_{\omega}=0$.
\end{proposition}
\begin{proof}
As the center of $\pi_1(M)$ is finite, and the flux is multiplicative when we replace a loop by a 
multiple, we may assume that the evaluation of a loop $\varphi_t$ in $\Symp_0$ whose flux we 
want to test bounds a $2$-disk $D$ in $M$. If $\gamma\subset M$ is any closed loop, the degree
$2$ homology class of $\varphi_t(\gamma)$ can be represented by a $2$-sphere $S$ obtained 
by surgering the torus along a meridian, using two copies of $D$. Now $\omega\vert_{\pi_2(M)}=0$
implies that 
$$
\langle \Flux_{\omega} (\varphi_t), [\gamma]\rangle = \int_{\varphi_t(\gamma)}\omega =\int_S\omega=0 \ .
$$
\end{proof}

\begin{example}
    Consider a surface bundle $X$ over a surface $B$, such that both the base $B$ and the 
    fiber $F$ have genus $\geq 2$. Then the second Chern number is nonzero, and the center 
    of $\pi_1(X)$ is trivial, so that $\Gamma_{2}=0$ by Theorem~\ref{t:kfl}. On the one hand, there are 
    many such $X$ which cannot satisfy the hard Lefschetz property (for any $\omega$), so that 
    we cannot conclude the vanishing of $\Gamma_{\omega}$ from Theorem~\ref{t:kfl}. On the other hand, 
    it is always possible to choose $\omega$ in such a way that it represents a bounded cohomology 
    class\footnote{This is implicit in~\cite{HK}, where it was used to show that $X$ always has positive 
    simplicial volume.}, in which case Corollary~\ref{t:sb} implies the vanishing of the flux group.
    
    Again this last argument does not cover all cases, because for suitable surface bundles $X$ 
    one can also choose $\omega$ so that its cohomology class is not bounded. Nevertheless,
    Proposition~\ref{p:P} always applies, because $X$ is aspherical and its fundamental group has 
    trivial center.
    \end{example}
    
\section{Final comments and remarks}

\subsection{Does the volume flux group depend on the smooth structure?}\label{ss}
   So far we do not know whether the non-triviality of the volume flux group depends 
   on the smooth structure, or not. Examples~\ref{ex:7D} and~\ref{ex:tori} are the 
   closest we have come to seeing a dependence on the smooth structure, but these 
   example are not conclusive. Most of the information we have derived from the non-triviality 
   of the volume flux group, about the fundamental group, homotopical irreducibility, 
   simplicial volume, and about the minimal volume entropy, is homotopy invariant. However, 
   this is not known for the minimal topological entropy, and is definitely false for the minimal volume. 
   Theorem~2 of~\cite{entropies} shows that in dimension $4$ there exist homeomorphic manifolds 
   such that one has vanishing minimal volume and the other one does not. It is also known that 
   the minimal volume depends on the smooth structure in higher dimensions.
   
\subsection{Remarks on Gottlieb groups}
Recall that the Gottlieb group $G(M)$ of a manifold $M$ is the image of the
evaluation homomorphism $ev_*\colon\pi_1(M^M)\longrightarrow\pi_1(M)$,
see~\cite{G0} and~\cite{LO,Oprea}. It will be clear to the experts that some 
of our arguments concerning volume flux groups depend only on the fact
that a loop of diffeomorphisms having non-trivial volume flux gives an element
of infinite order in the Gottlieb group. Indeed, the non-triviality of the Gottlieb group is 
enough to conclude that the Euler characteristic of $M$ vanishes, see~\cite{G0}, and for 
the irreducibility conclusion of Theorem~\ref{t:irred}. It is not clear whether 
the other consequences of a non-trivial volume flux follow from Gottlieb theory alone.
If $M$ has a Gottlieb element whose image under the Hurewicz map has infinite 
order in homology, then the simplicial volume of $M$ vanishes. However our proof 
of the vanishing of the minimal volume entropy certainly does not apply in this generality.
 
Once again the situation is better for $3$-manifolds. If a closed $3$-manifold without any 
fake cells has non-trivial Gottlieb group, then it is Seifert fibered, and up to multiples the
elements of the Gottlieb group are represented by circle actions. This follows 
from the Seifert fiber space conjecture, the final cases of which were settled by 
Casson--Jungreis~\cite{CJ} and Gabai~\cite{Gabai} independently. Our proof of 
Theorem~\ref{t:Seifert} did not need these deep results, because the existence of 
non-trivial volume flux implies that the manifold is Haken.

\subsection{Further developments}
Extended flux homomorphisms arose first in~\cite{KM} for the case of monotone symplectic
forms. There, a vanishing theorem for flux groups was proved as a byproduct of 
the search for extended flux homomorphisms, whereas in the present paper we obtain
many more vanishing theorems in the general situation, where an extended flux 
homomorphism may not necessarily exist. 

The results of~\cite{KM,KM2} illustrate how extended flux homomorphisms can help in understanding 
the homology of the groups $\Diff^{\alpha}$ as discrete groups. For many of the situations 
where we have proved the existence of extensions of the flux as a crossed homomorphism 
in this paper, one can try to imitate the constructions of~\cite{KM2} in particular in order to 
find new non-trivial cohomology classes on diffeomorphism groups made discrete. 

If an extended flux homomorphism $\widetilde{\Flux_{\alpha}}$ exists, then its kernel is an 
interesting subgroup of $\Diff^{\alpha}$ intersecting all connected components. In the special case when 
$\alpha$ is a symplectic form, this subgroup was recently studied by McDuff in~\cite{McDuffPrep}. She 
gives both an alternative proof of Theorem~\ref{t:mainA} for symplectic forms, and some elaborations on it.


\subsection*{Acknowledgements}
The first two authors would like to thank T.~Schick for a helpful discussion. We also acknowledge 
useful comments by B.~Hanke and by an anonymous referee.
This work was begun while the first author held a fellowship at the University of Munich funded by the 
{\sl European Differential Geometry Endeavour} (EDGE), Research Training Network HPRN-CT-2000-00101,
supported by The European Human Potential Programme. The paper was completed during visits of the second 
author to the University of Tokyo and to Stanford University. The third author is partially supported by JSPS Grant 
16204005. We are grateful to all these institutions for their support. 

\bibliographystyle{amsplain}

\end{document}